# PASSIVE SORTING OF ASTEROID MATERIAL USING SOLAR RADIATION PRESSURE

D. García Yárnoz,[*] J. P. Sánchez Cuartielles,[†] and C. R. McInnes [‡]

Understanding dust dynamics in the vicinity of asteroids is key for future science missions and, in the long-term, for asteroid exploitation. This paper analyzes the feasibility of manipulating asteroid material by means of solar radiation pressure. A novel method is proposed for passively sorting material as a function of its grain size or density, where solar radiation pressure is used as a passive in-situ 'mass spectrometer'. A simplified analysis shows that in principle this method allows an effective sorting of regolith material. This could have immediate applications for a sample return mission, and for industrial scale in-situ resource utilization to separate and concentrate regolith according to particle size or composition.

## INTRODUCTION

Asteroids are regarded as prime targets for space exploration missions. This interest is justified as asteroids are among the least evolved bodies in the Solar System and can provide a better understanding of its formation from the solar nebula [1]. Under NASA's flexible path plan [2], asteroids have also become attractive targets to be visited by crewed missions, with the benefit of not requiring the capability to land in and take-off from a deep gravity well. In addition, they may well be the most affordable source of in-situ resources to underpin future space exploration ventures [3].

---

[*] PhD Researcher, Advanced Space Concepts Laboratory, Dept. of Mechanical and Aerospace Engineering, University of Strathclyde, Glasgow G1 1XQ, UK.
[†] Research Fellow, Advanced Space Concepts Laboratory.
[‡] Director, Advanced Space Concepts Laboratory.

To date, in-situ observations of asteroids (e.g., Itokawa, Eros…) indicate that all Near-Earth Objects (NEOs) visited thus far, including very small bodies, are not bare lumps of rock [4]. A very fine layer of regolith material is likely to have a ubiquitous presence on most asteroid surfaces. The formation of this layer of regolith is usually explained by the effect of impact cratering and sandblasting through micro-meteoroid bombardment [5]. The presence of this fine dust, coupled with weak and irregular gravitational and electrostatic forces, increases the risk of triggering transient dust atmospheres during asteroid operations that can potentially degrade instrumentation, damage mechanisms and reduce visibility and communications. The future exploitation of asteroid material would need to take into account the dynamical behavior of dust under solar radiation pressure (SRP) in order to minimize the risk of such transient dust atmospheres. Considerable efforts have been made to understand the perturbing forces and space environment in the vicinity of cometary and asteroid bodies [6, 7]. These perturbing forces need to be considered and will have direct implications for the operations of spacecraft around and on small bodies. On the other hand, they also represent an opportunity, if engineered for practical benefit, to devise new types of highly non-Keplerian trajectories and novel methods for asteroid resource exploitation.

Extra-terrestrial resource exploitation is by no means a new idea. It was first proposed well over a century ago by the first pioneers of astronautics [8], and in the past decades it was given a comprehensive treatment by Lewis [9]. Ross further discusses the feasibility of extra-terrestrial mining applied to the NEO population [10]. The concept is presently back in the spotlight due to the founding of two companies with the final objective of mining asteroids: Planetary Resources[*] and Deep Space Industries[†]. If in-situ industrial scale exploitation is ever considered, various separation and material processing techniques would need to be implemented.

On Earth, industrial separation processes for mineral processing [11] range from the more traditional gravity concentration devices to numerous 'modern' methods including magnetic and electrostatic separation or, more recently, automated ore sorting [12]. Methods based on gravity or centrifugal separation are still used extensively in mineral processing as a first step to generate mineral concentrates for further treatment, or to discard waste, due to their simplicity and high capacity [13]. However, most of these gravity-driven separation methods are clearly no longer applicable in microgravity, or have a reduced performance, whereas others may require large-scale in-situ machinery.

There is an abundant literature on proposals for the exploitation and processing of lunar regolith [9, 14]. Magnetic separation techniques have been tested and proven useful on lunar simulants generated in the laboratory [15, 16]. Both reports show that paramagnetic pyroxene silicates and non-magnetic plagioclases can be effectively sorted with magnetic separation. So-called dry methods are effective down to particle sizes of 150 μm, whereas for smaller particles the use of slurries is needed [15], possibly due to cohesion. Further tests carried out in simulated lunar gravity on parabolic flights on tribocharged lunar silicate simulants [17] demonstrated the effectiveness of magnetic and electrostatic techniques in low gravity. Equivalent processes could be applied to asteroids [18], given the similarities in the silicate minerals present in both types of regolith and the vacuum and low-gravity environment. The higher ferro-metallic content in asteroid regolith would suggest that techniques based on magnetic separation are even more suited for asteroid resource exploitation to separate metals. Although this may be true, they require large, complex machinery for the separation, and for the previous steps of grinding and feeding. These methods would benefit from a prior regolith size separation or mineral concentration process. With this intention, new methods that take advantage of the low-gravity and vacuum environment of asteroids could be utilized.

---

[*] http://www.planetaryresources.com/
[†] http://deepspaceindustries.com/

This paper proposes and performs a feasibility analysis of one such novel method for sorting asteroid material, exploiting the dynamical interaction of regolith particles with solar radiation pressure. Separation is achieved by differential solar radiation pressure on ejected particles of different area-to-mass ratio. The concept is analogous to the separation process of 'winnowing' in agriculture, used for many thousands of years for separating grain from chaff due to differential atmospheric drag, again for materials with different area-to-mass ratio. This method has potentially attractive applications for large-scale industrial exploitation of asteroids, such as allowing a first, coarse, in-situ separation of different regolith particle sizes, or pre-concentration and separation of different materials based on their density. Future asteroid engineering and mining endeavors would benefit from this sorting technique, where solar radiation pressure is used as a passive in-situ 'mass spectrometer'. This process could be used in combination with, or as a first stage of a more complex process exploiting electrical or magnetic effects for more precise sorting.

**Exploiting Solar Radiation Pressure for Material Sorting**

Regarding the dynamical environment of asteroids, depending on an asteroid's size and its spin state, the effective ambient gravitational acceleration experienced by dust grains on small bodies can range from micro-gravity to milli-gravity [6], much lower than on the Moon. Under such conditions, the SRP perturbation becomes the largest non-gravitational force affecting single grains that have been lifted from the asteroid's surface, either naturally by micrometeoroid impacts or electrostatic forces, or artificially by mechanical means. Dust grains with a large area-to-mass ratio can escape from the asteroid [19], whereas those with smaller area-to-mass ratios will remain bounded. Their trajectories will nevertheless be significantly perturbed, even when ejected at low initial velocities. Based on this effect of differential SRP influence on dust grains, a method can be designed for passively sorting asteroid material as a function of grain size or density. The proposed idea consists of one surface element that collects and scoops loose regolith directly from the asteroid surface and expels it at a small velocity (either before or after grinding it), and one or several collectors that capture sorted particles as they fall back to the surface. A similar natural process, caused by electrostatic levitation [20], is believed to be responsible for the movement and concentration of fines (particles smaller than 40 μm) on shaded or shallow areas on asteroid surfaces. This process has been suggested as a possible explanation for the dust ponds observed at Eros [21].

In the following sections, simplified equations will be presented describing the trajectories of dust particles in the vicinity of an asteroid considering the third body perturbation of the Sun and the solar radiation pressure perturbation. In order to achieve a better understanding of the different regimes experienced by orbiting dust, a simple analytical formulation is applied to the problem. Both the analytical approximation and subsequent numerical propagations prove useful in studying the behavior of the dust particles. Preliminary conclusions can then be drawn from the analysis regarding the prospect of actively engineering and exploiting the forces experienced by dust grains in the vicinity of asteroids. Two separation strategies are presented, and the effect of uncertainties in the initial conditions on the separation is analyzed.

The preliminary analysis presented in this paper does not consider the modeling of additional perturbations, among them inter-particle forces, leaving it for future work. However, given the low gravity and vacuum conditions around asteroids, these additional forces acting on dust grains, particularly cohesion between individual particles, are likely to reduce the efficiency of the pro-posed sorting method. The implications of inter-particle forces are discussed by the authors in the final section of the paper.

**THE PHOTO-GRAVITATIONAL RESTRICTED 3-BODY PROBLEM**

The problem to be tackled can be modeled, in a first approximation, by the well-known photo-gravitational circular restricted three-body problem [22-24] applied to a spherical asteroid. The solar radiation pressure force acting on a particle is modeled with the standard cannon-ball approach:

$$\bar{F}_{SRP} = \frac{LQA}{4\pi c} \frac{\bar{r}_{SUN}}{r_{SUN}^3} \qquad (1)$$

where $L$ is the solar luminosity, $Q$ the solar radiation pressure coefficient, which depends on the material properties, $A$ the cross-sectional area of the particle, $c$ is the speed of light and $\bar{r}_{SUN}$ is the radio-vector from the Sun to the particle. The solar radiation pressure coefficient is 1 for a perfectly absorbing particle, and is equal to 2 for the case of ideal specular reflection. Unless otherwise stated, for the analysis in this paper the conservative value of $Q=1$ is assumed.

Assuming spherical particles of constant density $\rho$, the ratio of SRP perturbation with respect to the gravitational attraction of the Sun can be represented with the particle lightness number β given by:

$$\beta = \frac{LQ}{4\pi c \mu_S} \frac{S}{m} = \frac{3LQ}{16\pi c \mu_S r \rho} \qquad (2)$$

where $\mu_S$ is the gravitational constant of the Sun, $S/m$ the area to mass ratio of the particle, and $r$ the equivalent particle radius. Note that this lightness number has been defined with respect to the gravitational attraction of the Sun, and not the asteroid gravity as it is usual in literature. This parameter will prove useful to describe the different orbiting regimes of particles in the vicinity of an asteroid. Clearly, β is proportional to the particle's area-to-mass ratio, so for a fixed density, β increases when the particle radius decreases. Therefore, for small dust grains, SRP can provide a significant perturbing force.

For simplicity, it is assumed that the asteroid follows a circular orbit of heliocentric distance $d$ around the Sun. In a co-rotating frame with the origin at the barycenter of the system and with the x-axis pointing towards the asteroid (see Figure 1), the motion of a particle can be described by the following set of differential equations:

$$\begin{cases} \ddot{x} - 2\Omega_R \dot{y} = \Omega_R^2 \left( x - \frac{(1-\mu)(1-\beta)(x+\mu)}{\left((x+\mu)^2 + y^2 + z^2\right)^{3/2}} - \frac{\mu(x+\mu-1)}{\left((x+\mu-1)^2 + y^2 + z^2\right)^{3/2}} \right) \\ \ddot{y} + 2\Omega_R \dot{x} = \Omega_R^2 \left( y - \frac{(1-\mu)(1-\beta)y}{\left((x+\mu)^2 + y^2 + z^2\right)^{3/2}} - \frac{\mu y}{\left((x+\mu-1)^2 + y^2 + z^2\right)^{3/2}} \right) \\ \ddot{z} = \Omega_R^2 \left( -\frac{(1-\mu)(1-\beta)z}{\left((x+\mu)^2 + y^2 + z^2\right)^{3/2}} - \frac{\mu z}{\left((x+\mu-1)^2 + y^2 + z^2\right)^{3/2}} \right) \end{cases} \quad (3)$$

$$\mu = \frac{\mu_A}{\mu_A + \mu_S}; \quad \Omega_R = \sqrt{\frac{\mu_A + \mu_S}{d^3}} \quad (4)$$

where distances have been normalized with respect to $d$, $\mu_S$ and $\mu_A$ are the gravitational parameters of the Sun and the asteroid respectively, and $\Omega_R$ is the frequency of rotation of the two bodies (and the frame) around the barycenter. Higher order gravitational perturbations of the asteroid, which can be of great importance for irregularly shaped asteroids, are not taken into account at this stage of the investigation. Similarly, any inter-particle forces are neglected, and eclipses have also been ignored.

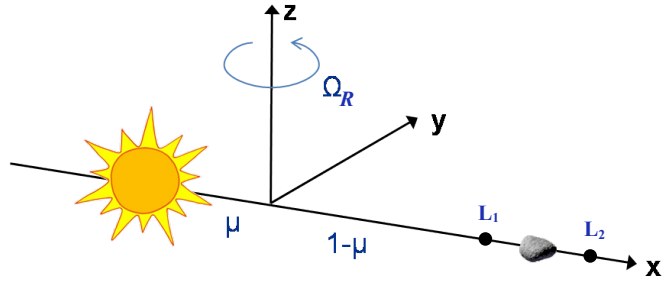

Figure 1. Schematic representation of the co-rotating frame with the origin at the barycenter of the Sun-asteroid system.

This system of equations has an integral of motion $C$:

$$C = -2E = 2U - 2T \tag{5}$$

where $U$ and $T$ are the potential and kinetic energy, which can be expressed as:

$$U = \frac{1}{2}\Omega_R^2 (x^2 + y^2) - \frac{(1-\mu)(1-\beta)}{\left((x+\mu)^2 + y^2 + z^2\right)^{1/2}} - \frac{\mu}{\left((x+\mu-1)^2 + y^2 + z^2\right)^{1/2}}$$

$$T = \frac{1}{2}(\dot{x}^2 + \dot{y}^2 + \dot{z}^2) \tag{6}$$

For particles ejected from the surface, it is possible to calculate zero velocity curves (corresponding to $C = 2U$) that depend on the lightness number $\beta$. Figure 2 represents the zero velocity curves for different particle sizes ejected with a fixed velocity from a hypothetical 10 km radius asteroid at 1 AU from the Sun, along with a set of example trajectories for a particular ejection site on the equator. The asteroid is assumed to be rotating with a 4 hour period around the $z$-axis, and the ejection velocity direction is selected radially outwards, normal to the asteroid surface. The average NEO density of 2.6 g/cm$^3$ [25] is considered for the asteroid, whereas estimates on the particle radius are provided assuming spherical grains of constant density of 3.2 g/cm$^3$, representing a relatively low density olivine. If the composition and structure of the asteroid is uniform, this implies a macro-porosity of 19%, close to the average S-type asteroid [26]. For ejection velocities above 11.2 m/s the zero velocity curves are open around L$_2$ for all values of $\beta$. An ejection velocity of 10.34 m/s was selected so that all particles with $\beta$ lower than 5x10$^{-3}$ (which corresponds to particles larger than approximately 35 µm radius) have closed zero velocity curves. This velocity increases with the rotation period (e.g. for a fast rotator with shorter 3 hour period, a correspondingly lower ejection velocity of 9.5 m/s has similar effects).

The range of lightness number selected covers from large boulders down to particles of size tens of µm. Smaller particles are likely to levitate naturally, or potentially escape when ejected due to the SRP perturbation (an example of this can be observed in Figure 2b). Two values of $\beta$ are particularly relevant: the $\beta$ for which the zero velocity curves open ($5\times10^{-3}$ in this particular case), which sets an upper bound in particle size for dust to escape, and the value of $\beta$ that ensures a re-impact before one revolution.

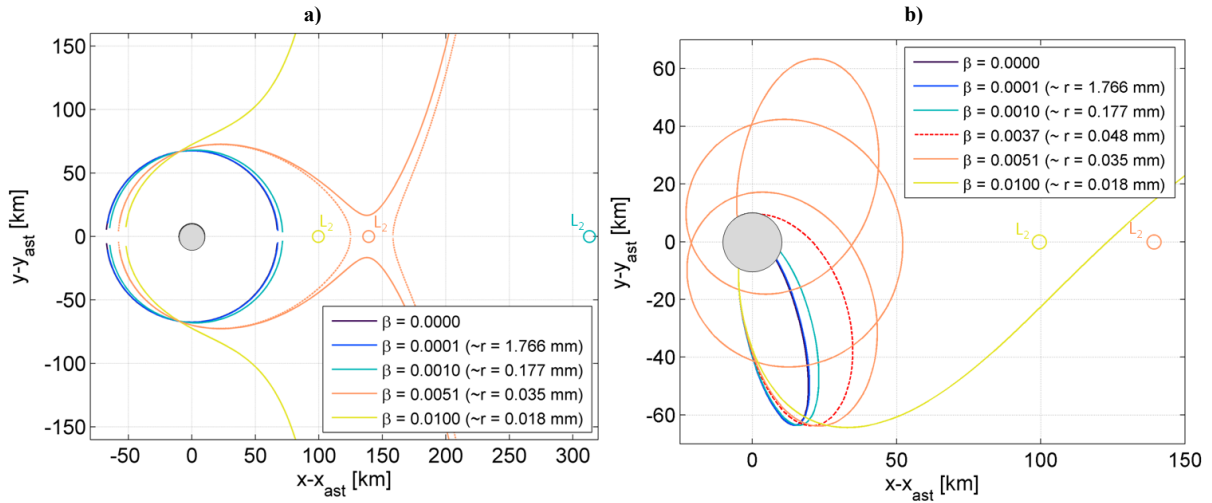

**Figure 2. Zero velocity curves (a) and trajectories in the co-rotating frame (b) for ejection velocities of 10.34 m/s from a 10 km asteroid with a 4 hour rotational period, for different values of $\beta$. Particle size estimation are given assuming spherical grains of constant density and $Q$=1. Zero velocity curves are open for all particles with $\beta$ larger than 0.0051. All particles with $\beta$ lower than 0.0037 re-impact before one revolution for this particular ejection site (dashed red line).**

The first value sets a theoretical limit for particles to escape, based solely on energy considerations. However, this provides little or no information in order to predict whether a particle would escape or not, as reaching escape conditions depends not only on the energy level, but also the solar longitude of the ejection site and the orbital geometry in general.

The positions in the *x*-axis of the two collinear libration points of the Sun-asteroid system, $x_1$ and $x_2$, can be calculated solving the system in Eq. (7) for each particular lightness number, and from there the associated value of the integral of motion *C* is obtained.

$$x_1 - \frac{(1-\mu)(1-\beta)}{(x_1+\mu)^2} + \frac{\mu}{(x_1+\mu-1)^2} = 0; \quad \rightarrow \quad C_1 = 2U(x_1)$$
$$x_2 - \frac{(1-\mu)(1-\beta)}{(x_2+\mu)^2} - \frac{\mu}{(x_2+\mu-1)^2} = 0; \quad \rightarrow \quad C_2 = 2U(x_2)$$
(7)

Given these values of *C*, and equating them to Eq. (5), substituting the potential energy of points on the surface of the asteroid, the ejection velocity required for the zero velocity curves to be open through the $L_1$ or $L_2$ points can be obtained. Because of the influence of the SRP perturbation, it is the $L_2$ point that offers the lowest energy for the zero velocity curves to be open, and thus also the lowest ejection velocity for escape trajectories of the dust grains. The main effect of SRP on the location of the libration points is to displace both the $L_1$ and $L_2$ points towards the Sun, resulting in an $L_2$ point closer to the surface of the asteroid (see Figure 2). Theoretically there is a particle size at which the $L_2$ point would be on the surface of the asteroid, and any smaller particles lifted from the surface with an infinitesimally small velocity at the correct time in a rotational period of the asteroid may escape. In reality, that region is in eclipse and the SRP would only affect particles that are still orbiting when the Sun comes into view.

It is then possible to calculate a guaranteed return velocity, for a given value of *β*, that ensures closed zero velocity curves and eventual re-impact with the surface. Figure 3 plots this velocity for a 10 km asteroid rotating along the *z*-axis assuming vertical ejection along the equator, where the asteroid's rotation provides the largest contribution to the total kinetic energy, and thus an easier escape at lower velocities. This guaranteed return velocity shows small variations with the longitude of the ejection site.

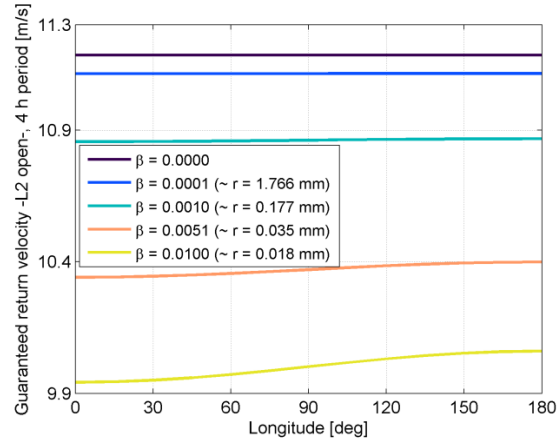

**Figure 3. Guaranteed return velocity for dust particles ejected radially outwards along the equator of a 10 km asteroid with a 4 hour rotational period, for different values of *β*. The longitude of the ejection site is measured along the equator with respect to the antisolar direction.**

In contrast, the value of *β* that ensures re-impact before one revolution is highly dependent on the point of ejection, and the relative geometry of the orbit with respect to the Sun. Certain ejection sites will have all particles directly re-impact under the perturbation of solar radiation pressure, whereas others will have a limiting *β* that allows multiple revolutions. Figure 4 represents the re-impact time in a longitude-latitude grid of the ejection site for the selected vertical ejection velocity of 10.34 m/s. The lightness number *β* is set to 0.0045 to ensure closed zero velocity curves and re-impact of all particles. Figure 4b gives the re-impact time in number of periods calculated with the initial osculating semi-major axis at ejection. The semi-major axis, and therefore also the period, is larger at the equator. For most of the surface of the asteroid the re-impact time is less than one initial period, while there is a region in the longitude's third quadrant that allows multiple revolutions. As expected, it is the region near the equator that has the highest probability of generating ejecta that perform more than one revolution. Also, the longitude's first and second quadrant, where the SRP acts against the velocity reducing the pericenter height, have shorter re-impact times that the third and fourth quadrant where the effect of the SRP perturbation contributes to raising the pericenter height.

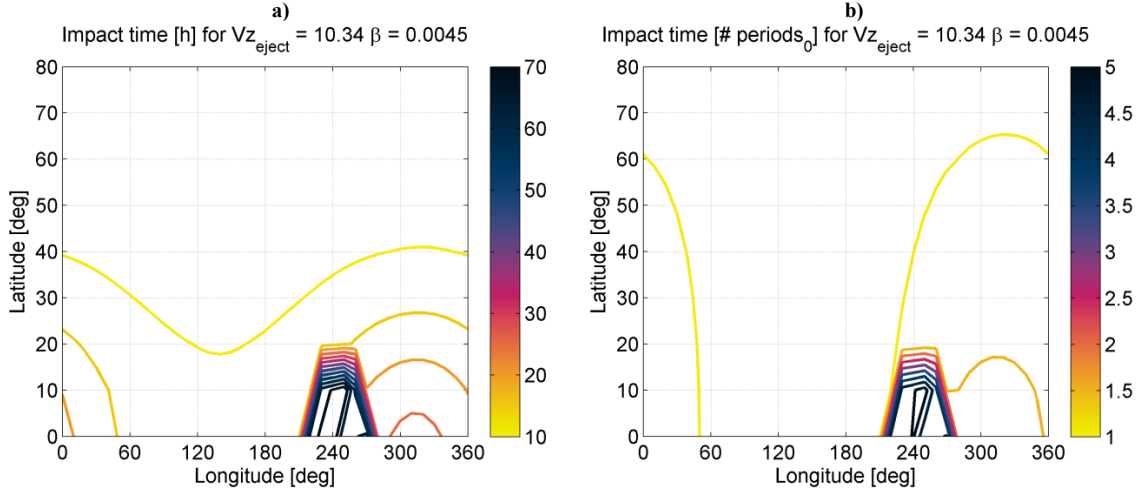

**Figure 4. Re-impact time on a 10 km asteroid with a 4 hour rotational period, for dust particles of $\beta = 0.0045$ ejected radially outwards with 10.34 m/s ejection velocity. The re-impact time is given in hours (a) and number of initial periods (b).**

It is difficult to draw more significant conclusions from this general form of the equations of motion, or to predict the re-impact of the dust without full numerical propagation for each particular case, which is time consuming. For this reason, a semi-analytical approximation is used in this paper to study the behavior and various regimes of ejected dust in the vicinity of an asteroid.

**HAMILTONIAN APPROACH**

The graph in Figure 6 is obtained by plotting the previously generated trajectories in an eccentricity-$\phi$ space, with the solar phase angle $\phi$ given by the following relation (see Figure 5a):

$$\phi = \Omega + \arctan\left(\frac{\cos i \sin \omega}{\cos \omega}\right) - \lambda_{SUN} + \pi \qquad (8)$$

where $\Omega$ represents the right ascension of the ascending node of the dust particle orbit around the asteroid, $i$ and $\omega$ are the inclination and argument of the pericenter, and $\lambda_{SUN}$ is the solar longi-tude. Initial ejection eccentricities range from 0.965 for an equatorial ejection, up to 1 for a polar ejection. It can be easily appreciated that the eccentricity increases along the trajectory when $\phi$ is larger than 180 degrees, and decreases when $\phi$ is between zero and 180 degrees. Only a few trajectories with initial solar phase angle $\phi$ around 90 degrees and latitudes close to the equator perform multiple revolutions. As an illustration, the evolution of a particular multi-revolution equatorial trajectory is plotted in Figure 6. The osculating ellipses are represented at ejection (A), re-impact (D) and two intermediate points when crossing the eccentricity value of 0.6. For this particular example, the eccentricity drops down to 0.5 before increasing back to values close to 1 causing a re-impact, while the solar phase angle evolves from 120 (A) to close to 240 degrees (D), as the pericenter rotates.

This eccentricity-solar phase graph resembles the phase space in the work by Oyama et al. [27] for the limiting case with infinite SRP-gravity ratio. Their analysis was based on an approach proposed by Hamilton and Krikov [28] to study the behavior of circumplanetary dust in a planar case, which consisted in orbit-averaging Lagrange's planetary equations over one revolution. This method was later used by Oyama et al. and various authors to describe applications for high-area-to-mass ratio spacecraft for Earth geo-magnetic tail exploration [27, 29], passive de-orbiting, and heliotropic orbits applications [30]. It can be shown that it is also well-suited to describe the trajectories followed by particles around asteroids under certain assumptions.

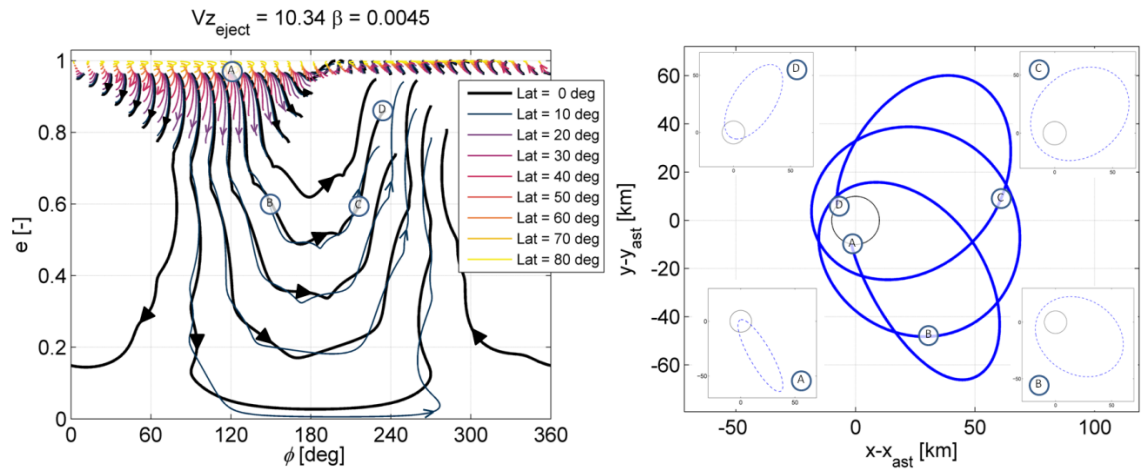

**Figure 5.** Eccentricity-$\phi$ plot of dust trajectories for different latitude-longitude ejection sites on a 10 km asteroid with a 4 hour rotational period, for dust particles of $\beta$ = 0.0045 ejected radially outwards with 10.34 m/s ejection velocity. An example equatorial trajectory is plotted, along with the osculating ellipses at ejection (A), re-impact (D) and two intermediate points along the trajectory.

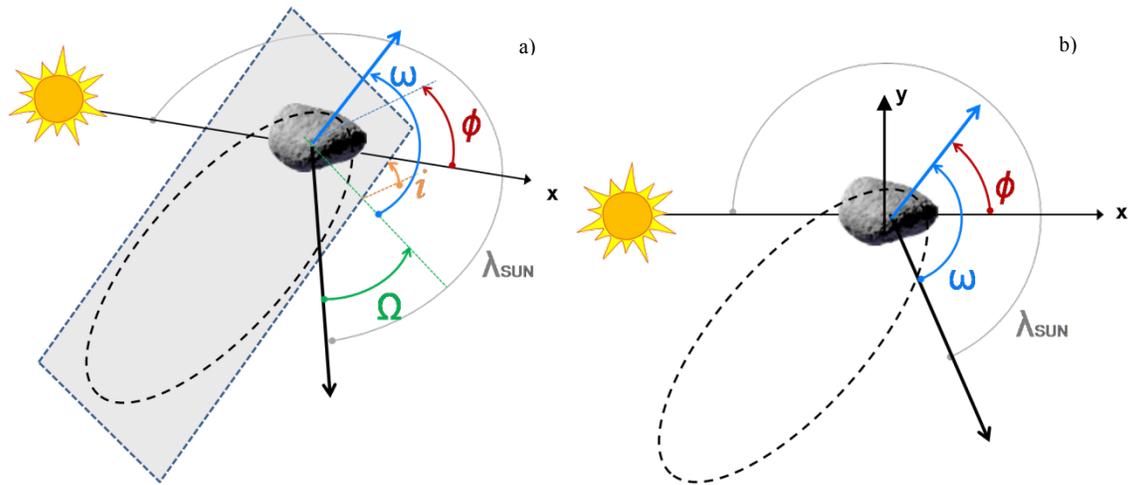

**Figure 6.** Definition of $\phi$ in the 3D case (a) and the planar case (b).

Following loosely Hamilton and Krikov's methodology, a planar case is assumed (see Figure 5b) in which the asteroid's rotational axis is perpendicular to the plane of movement of both the asteroid around the Sun and the dust particle around the asteroid. The phase angle $\phi$ between the anti-solar direction and the periapsis line is in the planar case simply given by $\phi = \omega - \lambda_{SUN} + \pi$. The dynamics of the dust under the influence of the solar radiation pressure perturbation and tidal forces caused by solar gravity can be described by the Hamiltonian [28]:

$$H = \sqrt{1-\bar{e}^2} + \frac{1}{2}A\bar{e}^2\left(1+5\cos(2\phi)\right) - C\bar{e}\cos\phi \tag{9}$$

where the coefficients $C$ and $A$ correspond to the SRP and tidal term respectively, and they can be expressed with the nomenclature followed by this paper as:

$$C = \frac{3}{2}\beta\sqrt{\frac{\mu_S}{\mu_A}\frac{\bar{a}}{d}} \qquad A = \frac{3}{4}\sqrt{\frac{\mu_S}{\mu_A}\frac{\bar{a}^3}{d^3}} \tag{10}$$

The eccentricity and semi-major axis that appear in Eq. (9) and Eq. (10) are orbit averaged values. This is an acceptable assumption in the case of circumplanetary dust, as the variation of the semi-major axis over one revolution is zero [28] and the eccentricity changes slowly. In the case of dust around asteroids, the excursions of the osculating semi-major axis from the mean and the variations in eccentricity in one revolution are much larger, introducing deviations from the analytical approximation, but the behavior of the system is still well described with the Hamiltonian approach. The evolution of the eccentricity and angle $\phi$ is then given by the following Hamiltonian system in non-canonical form, which uses the solar longitude as its independent variable.

$$\begin{cases} \dfrac{\partial \bar{e}}{\partial \lambda_{SUN}} = -\dfrac{\sqrt{1-\bar{e}^2}}{\bar{e}}\dfrac{\partial H}{\partial \phi} \\ \dfrac{\partial \phi}{\partial \lambda_{SUN}} = \dfrac{\sqrt{1-\bar{e}^2}}{\bar{e}}\dfrac{\partial H}{\partial \bar{e}} \end{cases} \tag{11}$$

While the initial argument of pericenter of the ejected dust can be arbitrarily chosen by selecting an ejection site/time, the rest of the initial osculating orbital elements $a_0$, $e_0$ and true anomaly $\nu_0$ of the ejected particles can be calculated as:

$$a_0 = \frac{R}{2 - \frac{R}{\mu_A}\left(v_{EJECT}^2 + \left(\frac{2\pi}{T_A}R\right)^2\right)}$$

$$e_0 = \sqrt{\left(v_{EJECT}\frac{2\pi}{T_A}R^2\right)^2 + \left(\frac{4\pi^2}{T_A^2}R^3 - \mu_A\right)^2} \bigg/ \mu_A \qquad (12)$$

$$v_0 = \arccos\left(\frac{\frac{4\pi^2}{T_A^2}R^3 - \mu_A}{e_0 \mu_A}\right)$$

with $R$ being the asteroid radius, $T_A$ its rotational period and $v_{EJECT}$ the ejection velocity, once again assumed normal to the surface of the asteroid.

Figure 7 plots the comparison between the numerical propagation of ejected dust trajectories, and the isolines of constant Hamiltonian that the orbit-averaged elements should follow in the analytical approximation. The ejection velocity of this particular plot is 9.5 m/s and the rotation period of the asteroid 3 hours. The initial eccentricity calculated with Eq. (12) is 0.93 for these ejection conditions. No equilibrium points are found in the phase-space, due to the high SRP perturbation when compared to the cases studied by Oyama et al. [27] and Colombo et al. [30]. This is consistent with the fact that no stable equatorial orbits can be found around small bodies when SRP is taken into account: the eccentricity starts eventually increasing up to values that cause a re-impact, or in the case of very high area-to-mass ratio, up to a hyperbolic escape. Still, the Hamiltonian approach correctly predicts the evolution of the eccentricity. For the initial conditions selected, there is an increase in eccentricity for all points with $\phi > 180°$, resulting in a decrease in pericenter height and immediate re-impact before one revolution of all particles ejected in two full quadrants of the asteroid.

The region around $\phi = 90°$ (which corresponds to ejection points with negative $y$ in the third quadrant of longitude as foreseen in Figure 4) contains both trajectories that re-impact and others that perform multiple revolutions. The eccentricity is decreasing in all cases, but only for a small range of initial phase angles does the pericenter height raise above the asteroid radius and the resulting trajectory performs multiple revolutions. To better illustrate this multi-revolution regime, a thick blue dashed horizontal line is plotted in Figure 7, which represents a critical eccentricity given by:

$$e_{CRIT} = 1 - \frac{R}{a_0} \tag{13}$$

corresponding to a value of 0.708 for the initial conditions selected in the figure. For eccentricities above this value and the initial semi-major axis, the osculating pericenter is below the asteroid surface, which results in a re-impact. As the variation of $a$ over one revolution is zero, and ejection/re-impact takes place close to pericenter, this approximation is accurate enough.

The apocenters and pericenters in the numerical trajectories have been indicated with X and O markers. Any pericenter (O marker) taking place above the critical eccentricity line implies a re-impact. It can be observed that whenever the eccentricity is lower than the critical value by the first pericenter, the dust particle manages to perform multiple revolutions, until after a few loops the eccentricity grows again above the critical one. The maximum reduction in eccentricity over one revolution does not take place though exactly at $\phi = 90°$, and the region of multi-revolution ejecta is thus not centered around it but shifted to the right on the plot.

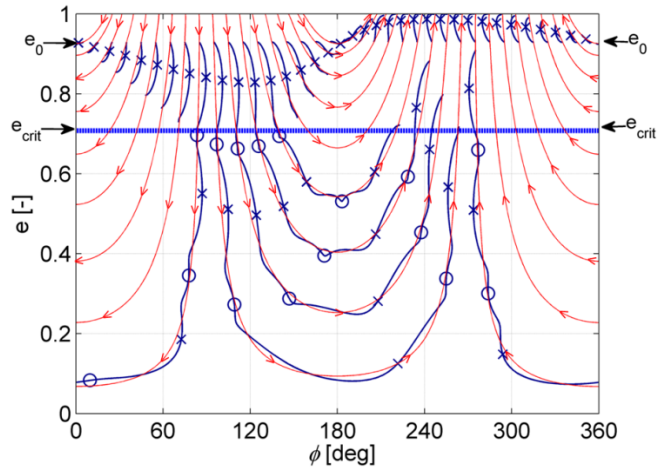

**Figure 7. Isolines of constant Hamiltonian (red with arrows) and numerical propagation of ejected dust (dark blue with X and O markers) plotted in the eccentricity-$\phi$ phase space. Apocenters and pericenters of the numerical trajectories are indicated with X and O respectively. Vertical ejection velocity is 9.5 m/s, $\beta$ = 0.0045, and the rotation period of the asteroid is 3 hours.**

A series of operational guidelines can then already be drawn from these results with regards to the selection of the extraction site, to determine the solar longitudes when operations are safer to avoid re-impact of dust on crewed missions or equipment. If transient dust atmospheres are to be avoided, solar longitudes close to $\phi = 270°\,(= -90°)$ would be preferred. Operations at other solar longitudes are still feasible if the forces used ensure the ejection velocity of dust stays well below the limit that allows multi-revolution trajectories.

**Time integration**

The Hamiltonian system in Eq. (11) can be transformed in full canonical form with the change of variable [28]:

$$k = \sqrt{1-\overline{e}^2} \tag{14}$$

resulting in:

$$\begin{cases} \dfrac{\partial k}{\partial \lambda_{SUN}} = \dfrac{\partial H}{\partial \phi} \\ \dfrac{\partial \phi}{\partial \lambda_{SUN}} = -\dfrac{\partial H}{\partial k} \end{cases} \quad (15)$$

It can be demonstrated that the tidal term $A$ is of the same order as the SRP only for distances of the order of the 20 asteroid radii, while the trajectories of interest stay bounded well below this distance. If the tidal term $A$ is neglected, Oyama et al. [27] showed that it is possible to integrate the system to obtain the time needed for a particle to travel along a line of constant Hamiltonian $H^*$ between two values of $k$ (or eccentricity), obtaining:

$$t - t_0 = \dfrac{-1}{\sqrt{\dfrac{d^3}{\mu_s}(1+C^2)}} \left( \arccos \dfrac{(1+C^2)k - H^*}{C\sqrt{1+C^2} - H^*} - \arccos \dfrac{(1+C^2)k_0 - H^*}{C\sqrt{1+C^2} - H^*} \right) \quad (16)$$

As a result, it is possible to plot isolines of transfer time on the phase-space graphs, in particular the isolines corresponding to the time until apocenter or pericenter, to determine if the next pericenter takes place before or after critical eccentricity is reached.

To obtain the orbital period it is also necessary to take into account the variation in argument of pericenter due to SRP. The Lagrange planetary equation for the derivative of time with respect to true anomaly, which includes a term related to the variation in $\omega$, is:

$$\dfrac{dt}{d\nu} = \sqrt{\dfrac{a^3}{\mu_A}} \dfrac{(1-e^2)^{3/2}}{(1+e\cos\nu)^2} \left[ 1 - \beta \dfrac{\mu_S}{\mu_A} \dfrac{a^2}{d^2} \dfrac{(1-e^2)^2}{e(1+e\cos\nu)^2} \left( \cos\phi + \dfrac{\sin(\nu)\sin(\phi+\nu)}{1+e\cos\nu} \right) \right] \quad (17)$$

Integrating over one revolution considering the mean values of semi-major axis and eccentricities, the orbital period can be approximated as:

$$T \approx 2\pi \sqrt{\dfrac{\bar{a}^3}{\mu_A}} \left( 1 - \dfrac{1}{2} \beta \dfrac{\mu_S}{\mu_A} \dfrac{\bar{a}^2}{d^2} \cos\phi \left( \dfrac{12 + 13\bar{e}^2}{4\bar{e}} \right) \right) \quad (18)$$

A similar approach can be followed to obtain the first semi-period until apocenter. Plotting on the phase space the isolines of Eq. (16) corresponding to the time from initial true anomaly to apocenter and pericenter calculated analytically (see Figure 8), they are a good match for the apsides points calculated with the numerical propagation. This allows an accurate analytical prediction of the conditions for a particle to perform multiple revolutions. The region in the phase space where the isolines of time to pericenter is below the critical eccentricity line corresponds to multi-revolution trajectories. It is important to note that the critical eccentricity varies with the semi-major axis along the orbit, so there can be pericenters close to the limiting values on either side that may impact or fly-over depending on the exact value of the semi-major axis at the re-impact time for the two extreme cases.

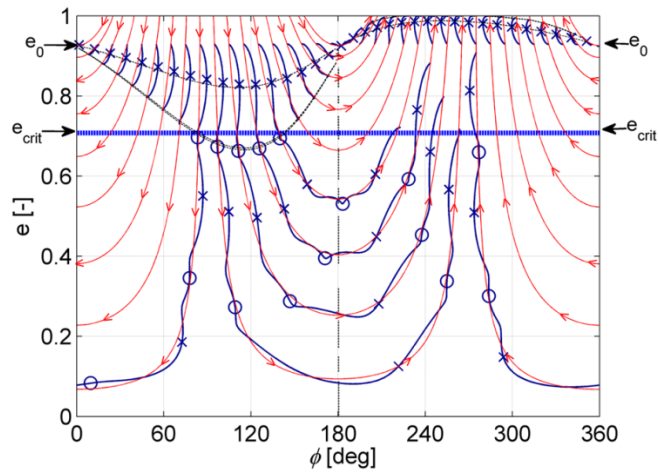

**Figure 8. Apocenter and pericenter analytical time estimation (dotted black lines) on the eccentricity-$\phi$ phase space. Vertical ejection velocity is 9.5 m/s, $\beta = 0.0045$, and the rotation period of the asteroid is 3 hours.**

**MATERIAL SORTING APPLICATIONS**

One of the benefits of the differential effect of solar radiation pressure on ejected dust particles is the possibility to engineer these forces in order to passively separate material as a function of $\beta$. This processing of material can be considered either for separation of the same material as a function of the grain size on an asteroid of uniform composition (the larger the grain, the lower the $\beta$), or alternatively, after a grinding process to reduce all materials to a similar grain size, as a method of separation of two materials with different density (again the higher the density of the material, the lower the $\beta$). The SRP sorting concept takes advantage of the low gravity on asteroids, which would render other techniques such as classical gravity concentration processes unfeasible.

Two possible strategies for collecting the separated material can be devised (see Figure 9): a hovering spacecraft that collects on-orbit the material that has been lifted from the surface by a surface element or rover, or several collection points on ground at pre-calculated distances. Each method has preferred ejection points on opposite sides of the asteroid ($\phi = \pm 90°$). The surface element can then collect material and, if needed, grind it for one asteroid revolution, and then eject it at the appropriate time depending on the selected strategy.

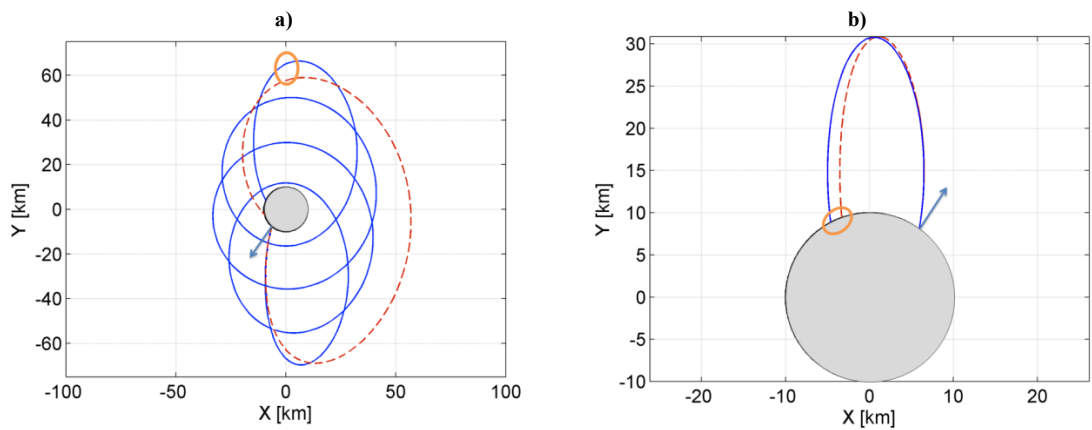

**Figure 9. Schematic representation of separation strategies with a hovering spacecraft collection point (a) and a ground-based collection point (b) for a heavier particle (solid blue) and a lighter one (dashed red). X-Y axes are parallel to the co-rotating frame axes and centered on the asteroid.**

**On-orbit collection**

In the first method, the preferred ejection point is close to $\phi = 90°$, where the SRP contributes to increasing the pericenter height and multiple revolutions are possible. The collection point should be hovering at a certain distance on the Y-axis. Once an ejection velocity has been selected, there is a minimum value of $\beta$ that will avoid re-impact in the first revolution, essentially discriminating a maximum size of the grains of interest (larger grains would fall back onto the asteroid). The rest of the material would travel following isolines of constant $H$ until eventually reaching the desired height on the Y-axis for collection. Figure 10 shows on an eccentricity-$\phi$ phase plot the lines that correspond to a particular height over the Y-axis on the asteroid assuming a constant mean semi-major axis. Only for certain ejection sites close to $\phi = 90°$ does the eccentricity evolution along the constant Hamiltonian isolines allow the particles to reach heights over 55 km on the Y-axis after a number of revolutions. The thick blue line shows a propagated trajectory of this type, where the eccentricity decreases down to 0.1 (almost circular), before it increases again when the phase angle shifts to values close to 270°. The sixth apocenter for this trajectory is well above 60 km.

The passive separation takes place in time, as larger particles (particles with smaller $\beta$) require more revolutions and thus longer times to reach the collection area. This is shown in Figure 11a, where the time to come back to the initial eccentricity levels increases with decreasing $\beta$. The plot shows the time evolution of the eccentricity for different values of $\beta$ and the same ejection site at the equator and $\phi = 90°$, a rotation period of 3 hours, and an ejection velocity of 9.5 m/s. Particles with $\beta \leq 0.004$ re-impact before the first pericenter.

It is possible to analytically estimate the time until collection by calculating with Eq. (16) the time to reach the eccentricity of intersection of the ejection point isoline with $\phi = 180°$ or $\phi = 0°$. Obtaining the intersection point is straightforward by substituting the phase angle $\phi$ in Eq. (9). Figure 11b shows a comparison between the analytical estimates and the numerical propagated trajectories. The collection point for each $\beta$ is assumed as the point when the eccentricity reaches the value at the first apocenter again, and is indicated with markers in Figure 11a. The first lighter particles arrive one day after ejection, while the heavier particles can take over 70 hours.

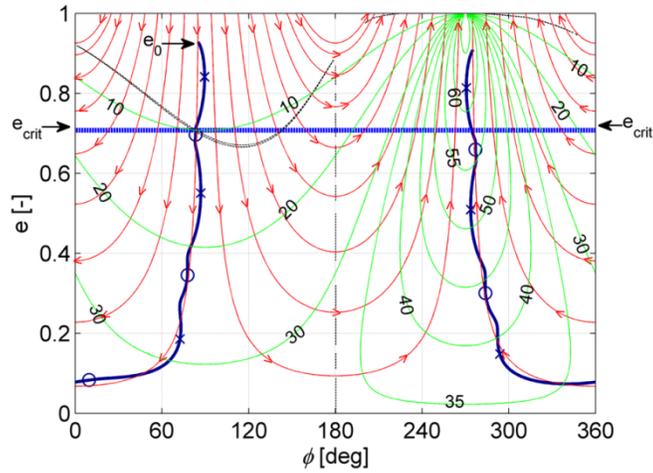

**Figure 10.** Eccentricity-$\phi$ phase space graph with isolines of constant height over the Y-axis in km (number-labeled green contours). Vertical ejection velocity is 9.5 m/s, $\beta$ = 0.0045, and the rotation period of the asteroid is 3 hours. For ejection sites close to $\phi$ =90° trajectories can reach heights over 55 km on the Y-axis when the phase angle shifts to 270° after a number of revolutions.

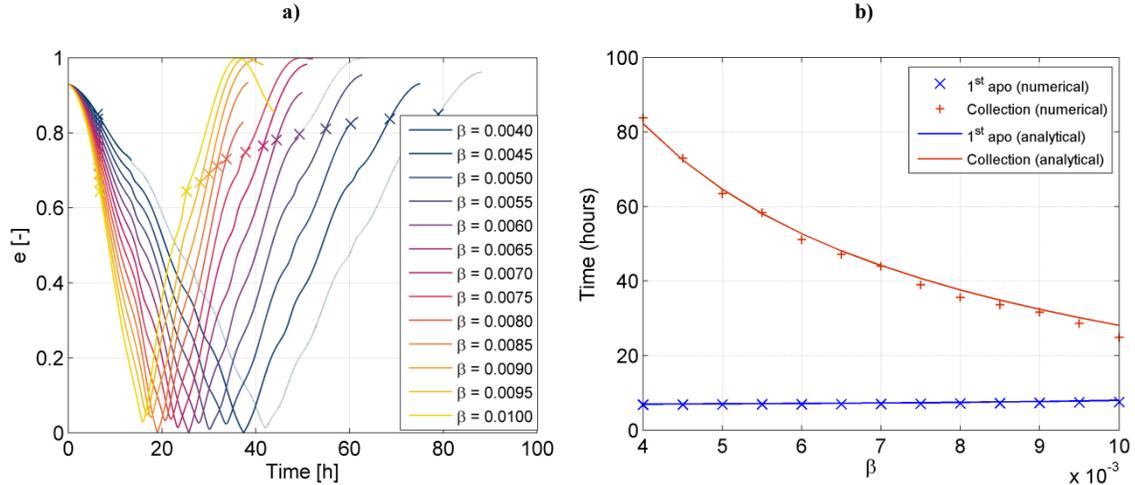

**Figure 11.** a) Evolution of the eccentricity with time for $\phi = 90°$ and different values of $\beta$. Markers indicate the first apocenter and collection point. Discontinuous lines indicate there has been a re-impact on the surface. b) Time to collection and analytical approximation as a function of $\beta$.

The main difficulty concerning this strategy is the large variation in height at the Y-axis crossing for orbiting particles of different $\beta$. Furthermore, simplifications in the model, such as the omission of the eclipse times and not including other perturbations such as higher order gravitational terms, make the predictions in time much less accurate and the dispersion at the collection point higher. A non-planar case would have the added difficulty of the evolution of the orbital plane. Very large collectors (of the order of hundreds of meters for an asteroid of radius 10 km) would be required in order to capture a significant number of particles, or an elaborate station keeping strategy that changes height with time would be needed in order to compensate for these variations. This greatly increases the complexity of operations for a hypothetical collection on-orbit. The ejection velocities involved are also higher than for a ground based collection point.

**On-ground collection**

The second method envisages a series of collector points spaced on the surface of the asteroid (or a collector band or strip extended over some distance). The differential separation of particles as a function of $\beta$ is performed in space, rather than in time. The preferred ejection point is in this case close to or equal to $\phi = 270°$ (pericenter on the negative Y-axis), which corresponds to trajectories where the SRP reduces the pericenter height and thus re-impact is assured before one revolution for most values of $\beta$.

Figure 12 shows the distance between re-impact points for two regolith particles of different sizes and densities as a function of the ejection velocity. The ejection velocity in the X-axis has been scaled with the radius of the asteroid. For asteroids in the size range of 100 m to 10 km the required ejection velocities for the same separation densities scale well with the radius. Only in the case of fast rotators (blue lines with 2.5 hour period) if the ejection velocity increases there is actually a clear bifurcation between the separation at a 1 km asteroid (dashed) and the 10 km one (continuous line) for the range of velocities plotted, and there can even be cases with particles escaping. For a smaller 100 m asteroid the bifurcation takes place at lower relative ejection velocities but in all cases it is well above the 1 m separation horizontal line.

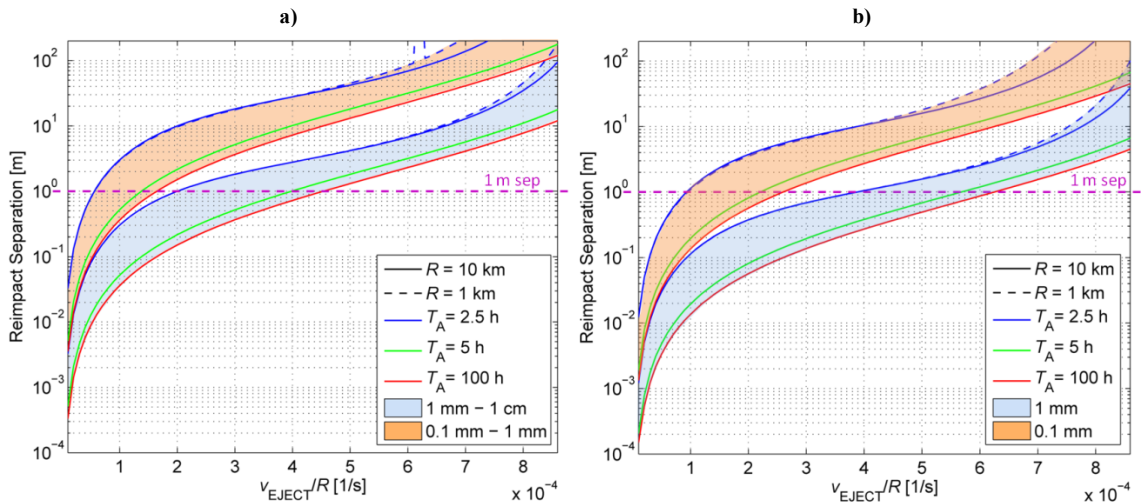

**Figure 12. Re-impact point separation between 1 mm-1 cm (blue) and 0.1 mm-1 mm (orange) grains of homogenous density of 3.2 g/cm³ (a), and between grains of different densities (2.68 and 3.74 g/cm³) of the same size, 1 mm in blue, 0.1 mm in orange (b). Particles are ejected with a phase angle $\phi$ of 270 degrees for different asteroid rotation periods. Dashed and solid lines indicate the 1 km and 10 km asteroid respectively. Only for the short period case (2.5 h) does the 1 km line deviate from the 10 km one for this range of ejection velocities.**

Figure 12a assumes spherical particles of constant average density of 3.2 g/cm³ (homogeneous asteroid of a low density olivine) for the suggested ejection site and different rotational periods of the asteroid. If the desired separation between two such particles of size 1 mm and 1 cm is 1 m, the graph shows that for an asteroid rotating with a 2.5 hour period the required ejection velocity would be 2 m/s on a 10 km asteroid. This velocity would more than double in the case of a non-rotating asteroid (~ 4.5 m/s). If the size of the particles desired for separation is one order of magnitude lower (0.1 mm and 1 mm), the ejection velocities range from 0.5 to 1.7 m/s for the same separation distance. These velocities are almost one order of magnitude lower than the ones suggested for a hovering spacecraft collection point.

Figure 12b assumes instead a differentiated asteroid with two materials to be separated. The regolith is assumed to have been previously ground to grains of the same size (1 mm or 0.1 mm).

Assuming 1 mm grains composed of plagioclase (average density of 2.68 g/cm$^3$) and a denser olivine or pyroxene (3.74 g/cm$^3$), the ejection velocities required for a 1 m separation on a 10 km asteroid range from 4 to 6.2 m/s depending on its rotation rate. For finely ground 0.1 mm particles, these velocities are reduced to 1 to 2.6 m/s. As the required ejection velocities scale with the radius, on a 1 km size asteroid they would be one order of magnitude lower.

In addition to the lower ejection velocities, other benefits of the ground base collection when compared to on-orbit collection is that eclipses have little or no influence in the trajectories of the particles as the preferred ejection sites result in effectively eclipse-free trajectories except for a short interval at ejection. Other perturbations, such as higher order harmonics of the gravity field for the usually irregularly shaped asteroids would affect all particles in a similar form, regardless of the $\beta$ value, so the differential effect of SRP would still cause a separation in re-impact points of the same order. Non-planar trajectories would be affected in a similar way over one revolution, and the only concern of an ejection point away from the equator would be a higher required ejection velocity, similar to the case of a slowly rotating asteroid.

**EFFECT OF UNCERTAINTIES**

One of the concerns that arise when evaluating this method is the sensibility of the separation to errors or uncertainties in the ejection conditions. As the relative difference in grain size or density is small, errors in the velocity at ejection may induce greater dispersion on the grains than the solar radiation pressure perturbation.

The separation caused by velocity errors was calculated for various asteroid sizes and rotation periods. Figure 13 shows the re-impact point separation for an error in the velocity modulus of 1% (a) and an error in the ejection direction of 0.33 degrees (b). The color patches indicating the separation by solar radiation pressure calculated in the previous section are superimposed for comparison. It can be observed that there is a limiting asteroid size for each ejection velocity for which the SRP induced separation is no longer effective compared to errors in velocity.

Focusing on specific examples on Figure 13a (indicated by vertical lines), for the case of a fast rotating asteroid with a period of 2.5 h (blue square ①), the required ejection velocity for a separation of 1 m is equal to 0.06 x $R$ m/s. With velocity modulus errors of 1%, the separation between equal particles would be of the order of 20 m for a 10 km asteroid, 2 m for a 1 km one, and 0.2 m for asteroids of 100 m radius. For a slowly rotating asteroid (period of 100 h, ②), the separations due to errors in velocity modulus are of the order of 0.5, 0.07 and 0.02 m for asteroids of 10, 1 and 0.1 km radius respectively. While the dispersions in a slow rotator for asteroids up to 1 km seem acceptable for implementing the SRP particle sorting, they would render it useless for larger asteroids of 10 km radius, or for fast-rotators.

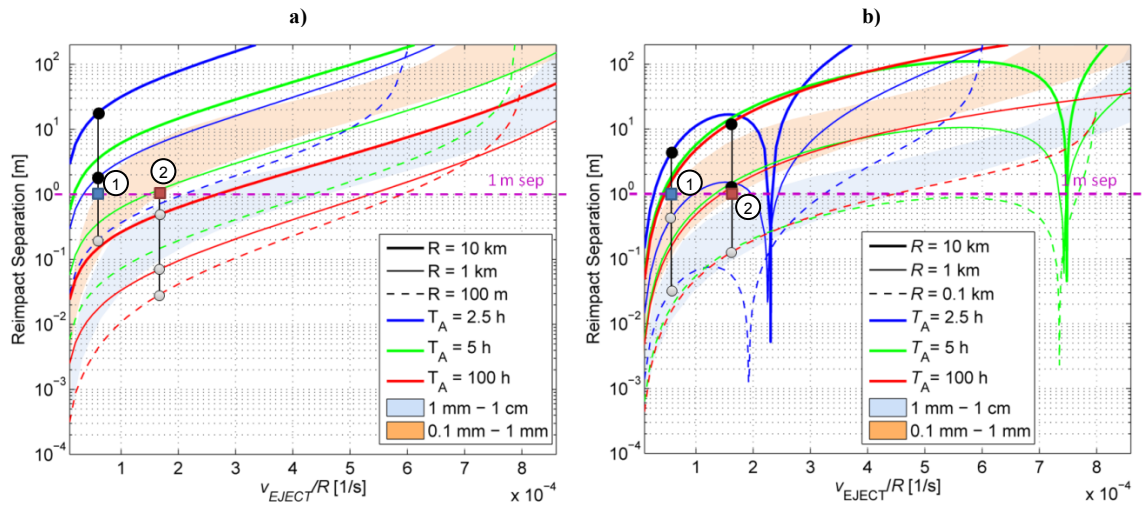

**Figure 13. Re-impact point separation for errors in ejection velocity modulus of 1% (a) and errors in ejection angle (in-plane) of 0.33 degrees (b). Line-styles indicate asteroid size; colors indicate its rotation period. The patches indicating the separation by differential SRP are also included for reference.**

Figure 13b shows a similar plot for errors in the velocity angle (only in-plane errors were considered). The separation induced by errors in angle has actually a minimum at a particular veloci-ty, which lies close to the required ejection velocity to obtain an orbital period for the grains equal to the rotational period of the asteroid. Close to this velocity, an error in the ejection angle causes an error in semi-major axis that barely modifies the re-impact time. As the particles perform a revolution at approximately the same time as the surface, there is a particular velocity for which errors in angle would result in the same re-impact point at the surface, only slightly before or after in time. When considering the same two examples as in the previous case, the separation induced by errors in ejection angle is of the order of 4, 0.4 and 0.03 m for the fast-rotator case ①, and 10, 1 and 0.1 m for the slow-rotator ②. This seems to indicate that only asteroids of up to 100 m radius (or slightly larger) are good candidates for the sorting method proposed.

In order to have a more comprehensive analysis of the effect of additional uncertainties, such as variations in the grain size or the densities of the materials, Monte Carlo simulations with 10000 shots were run for a 100 m radius fictitious asteroid rotating with a period of 5 h and with regolith that has previously been ground to a particle radius of 100 μm.

For the composition of the asteroid, an ordinary chondrite S-type asteroid is assumed, containing mostly silicates (in general olivine-pyroxene with densities ranging from 3.2 to 4.37 $g/cm^3$) and a few traces of metal (of density 7.3 to 7.7 $g/cm^3$) [26]. Table 1 lists the material composition as well as the assumed mean density and dispersion for 5 selected materials: Fe-Ni metallic grains, a high density orthopyroxene (opx), two olivine (ol) silicates of medium and low density and a plagioclase silicate. The average grain density is of 3.52 $g/cm^3$ for this particular mix, which is well within the range of ordinary chondrite meteorites of L or H type [31, 32].

Table 1. Regolith composition and densities. opx: orthopyroxene, ol: olivine

| Material | Percentage in regolith | Mean density [*], [26] (g/cm$^3$) | 1-σ Std. dev. (g/cm$^3$) |
|---|---|---|---|
| Fe-Ni | 2% | 7.50 | 0.07 |
| High density opx | 15% | 3.95 | 0.10 |
| Medium density ol | 50% | 3.50 | 0.10 |
| Low density ol | 28% | 3.20 | 0.05 |
| Plagioclase | 5% | 2.68 | 0.04 |

Table 2 presents additional variables with uncertainties in the Monte Carlo simulation. The velocity modulus is selected to obtain a nominal separation of 1 m for an asteroid with a rotation period of 5 hours for 100 μm particles of the materials selected in the previous section (see Figure 13b). The errors in velocity modulus follow a normal distribution with a 3-σ uncertainty of 3% of the velocity modulus. Errors in ejection angle are assumed in two orthogonal directions (along the longitude and latitude) with 1 deg 3-σ. The particles in the regolith are assumed to be ground prior to ejection down to a fine dust with radius of 100 μm. The distribution of particle size is assumed to follow a lognormal distribution of parameters $\mu_{log}$ and $\sigma_{log}$ of -9.21 and 0.05, which corresponds to a mode in particle radius of 100 μm. The mean value is slightly higher and the standard deviation is approximately 5 μm, as can be seen in Table 2.

Table 2. Uncertainties in Monte Carlo simulation.

|  | Mean | 1-σ Std. dev. |
|---|---|---|
| Velocity modulus (cm/s) | 2.35 | 1% |
| Error angle (deg) | 0.00 | 0.33 |
| Particle radius (μm) | 100.38 | 5.22 |

---

[*] http://webmineral.com last accessed 29/01/2013

Figure 14 shows the results of one Monte Carlo run with 10000 shots. Figure 14a represents the trajectory in the Sun-asteroid co-rotating frame, with the Sun in the negative X-axis. Metallic particles are the least affected by the SRP and fall closer to the trajectory without perturbations (red dashed line). In Figure 14c the same trajectories are plotted in a local horizontal frame with the $x_{loc}$-axis tangent to the longitude lines on the asteroid (in this case the equator) and the $z_{loc}$-axis normal to the surface. The Sun direction rotates but it is also close to the negative $x_{loc}$-axis due to the selected ejection point with phase angle $\phi = 270°$. From the point of view of the ejector (i.e., Figure 14c), the particles start travelling upwards along the local vertical up to a height of about 5 meters, and then start falling behind as the asteroid rotates, with the closest heavier particles falling around 3 m in the anti-solar direction and the lightest particle considered re-impacting 7 m away from ejection. A trajectory without SRP is also plotted (red dashed line) for comparison. The re-impact points are represented on Figure 14b and d. Figure 14b shows the displacement along the $x_{loc}$-axis for different densities, while Figure 14d plots the re-impact point in the local horizontal frame. The theoretical re-impact point without SRP is also indicated in both plots. The main effect of errors in velocity in the $y_{loc}$ direction is a displacement of the re-impact point in the same direction. Its influence in the separation or mixing of particles is limited.

The separation as a function of density is particularly effective for metallic particles, due to their well differentiated density. There is a much less clear separation among different silicate materials, although a gradient in density along the local *x*-axis is evident, from heavier pyroxene-olivine mixtures to lighter silicates of the plagioclase-feldspar family. The distances involved are also acceptable from an engineering point of view. A 10 meter collector band extended from the ejection module could be deployed to collect material, or a rover could sweep the re-impact area in strips perpendicular to the $x_{loc}$ direction.

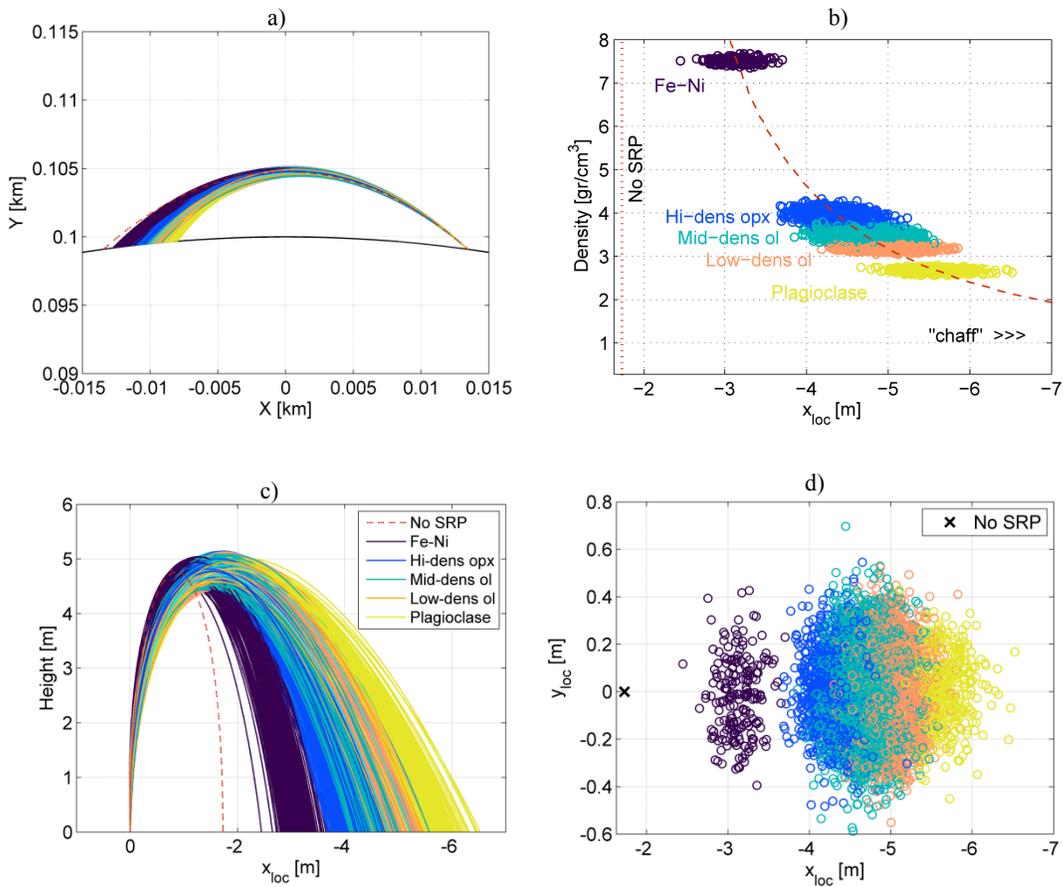

**Figure 14. Monte Carlo run with 10000 ejected particles of different materials. Plot a) shows trajectories in the Sun-asteroid co-rotating frame; b) shows the separation as a function of density; c) plots the height versus the local horizontal frame separation along the *x*-axis; d) shows the re-impact points in a local horizontal frame.**

**DISCUSSION AND MODELING LIMITATIONS**

The analysis presented so far represents a first feasibility study of the concept of passive sorting of material on asteroids by means of solar radiation pressure. This analysis already shows that ground based collection of material seems to be a promising technique. It would be particularly suited for concentrating materials with well differentiated densities, or to coarsely separate grain sizes.

In the first case of material separation as a function of density prior grinding is needed. A possible implementation would consist of regolith grinders that operate for a full asteroid revolution and eject the material only at the desired solar longitude. The collectors would be static, placed a few meters away or even attached to the ejector. Similar to mining terrestrial processes, the collectors could consist of simple mounds, strips of material, or "pools" on the surface where material accumulates for later collection or processing. The process can be repeated iteratively with the material in these mounds being ejected again from a different location for further refinement.

On the other hand, for regolith particle size separation, no grinding is required. The regolith collected could be directly ejected by a moving system that advances ideally at the same rate that the asteroid rotates (at a rate of a few millimeters or centimeters per second depending on the asteroid period). A first line of rovers could rake the regolith in advance, while a second line that maintains its position close to the terminator on the shaded region at the correct solar longitude scoops it and ejects it continuously. The collectors need also advance at the same rate and in coordination with the ejectors, maintaining their relative position in the lighted area after sunrise close to the terminator to receive the ejected material. However, surface mobility at the required speeds is at best questionable in the microgravity conditions and rough uneven surfaces of asteroids [33]. Legged rovers with gripping mechanisms are currently under development for movement in very steep surfaces on Mars or the microgravity conditions of asteroids [34, 35]. They would enable steady locomotion on asteroids and the ability to traverse obstacles, but unlikely at the required speeds to allow a continuous ejection system that keeps up with the asteroid rotation. Large industrial ventures may consider the installation of complex machinery and systems such as rails to enable locomotion on asteroids for exploitation.

Alternatively, both for density and regolith size separation, a possible solution is moving rovers or excavators that return to a designated ejection point at the appropriate time to eject the collected and/or ground material. These rovers could use the already mentioned gripping-motion systems, but at much lower speeds. The ejection sites could be preselected to take advantage of the local topography for better separation or ease of collection.

The main benefit of the separation method in itself is the simplicity of its physical principle. This is an advantage for implementation, which can be modular and scalable, with potential for high throughput. Because of this, it has potential applications for the concentration of materials or grain sizes prior to or in combination with more sophisticated processing techniques specific for each material type or particle size. It would comprise the first step of an industrial process similar to gravity concentration on Earth. Nevertheless, further work on the numerical validation of these semi-analytical results is desirable in order to ensure the efficacy of the method presented in a more complete dynamical model of a small body.

First of all, the use of the CR3BP for numerical propagation instead of a more complex model including the eccentricity of the asteroid orbit around the Sun has a minor influence on the results, as the timescales of the trajectories considered are short compared to the period of a NEO orbit around the Sun. In addition, the analytical results correspond to a simplified planar case and, although an extension to low inclination is possible and the general behavior is not expected to change, a three-dimensional model is required for large inclination trajectories for ejecta from higher latitudes or from asteroids where the rotation axis is not perpendicular to the Sun-line. For these cases, the closed-form solution of the radiation pressure approximation (RPA) by Richter and Keller [7, 36] could allow an analytical extension of the problem to 3D. The RPA approach has already been used in the literature to study stable orbits of ejecta among small bodies [36]. It may however not be best suited to study the trajectories of interest, as the model fails if the number of revolutions is small, as in this case.

The irregular shape of asteroids has also an important effect on the method proposed, not only due to unmodeled gravitational perturbations that may render some of the multi-revolution trajectories presented unfeasible [37, 38], but also due to the changes in the re-impact point position and direction as a function of the local geometry. Higher order gravitational terms should affect all particles independently of their lightness number, so their influence on the SRP induced separation is expected to be limited for the on-ground collection method where a full revolution does not take place. However, the shape of the asteroid will change the departure and re-impact conditions significantly.

Regarding the solar radiation pressure force modeling and the material properties, the same reflectance was assumed for all materials by setting the radiation pressure scattering coefficient $Q$=1, corresponding to complete absorption. The lightness number $\beta$ is proportional to $Q/\rho$ and material dependant variations in the radiation pressure coefficient would affect the separation between particles. However, reflectance spectra analysis from meteorite metal-silicate mixtures [39] show that the absolute reflectance of orthopyroxene and olivine silicates is greater than that of denser metallic mixtures, thus having the effect of a greater gap in $\beta$ values, and theoretically increasing the re-impact separation between metallic and silicate particles. The simple cannon-ball model used for SRP is considered accurate enough to study particle evolution, where the shape and effective area and the attitude of the particle are not known. In addition, as pointed out earlier, eclipses would not affect the on-ground collection scheme as preferred ejection points result on effectively eclipse-free trajectories (except for the ejection instant). However, if strategies involving multi-revolution trajectories are implemented, eclipses need to be considered.

They influence the evolution of eccentricity (and all other orbital elements) in the phase space plots and may result in earlier or later re-impact.

Inter-particle forces are of extreme relevance when separation methods are to be implemented, and they need to be taken into consideration. Scheeres [6] provides a scaling of these forces as a function of particle size that proves useful to discriminate between the various perturbations.

Self-gravity can in general be ignored for particle sizes below a few centimeters. That is not the case for electrostatic forces. It was theorized [20, 21] in fact that they provide a natural mechanism for segregating and accumulating particles smaller than 100 μm on craters and shaded areas through electrostatic levitation, if sufficient charging time is provided. These predicted smooth deposits of material in ponds were confirmed by observation at Eros [21], and the levitation and transport mechanism was reproduced by simulation [40, 41]. The proposed sorting method is in fact an enhanced extension of this natural 'electrostatic winnowing' for larger particles, with the ejection of particles not necessarily performed by electrostatic means.

The ejection sites selected in both cases are close to the terminator of the body, where electrostatic forces are expected to be larger [6], and dust levitation is most likely to take place naturally. The use of electrostatic forces as an aid or as the main ejection mechanism could therefore be considered, using artificial electric fields to accelerate tribocharged particles. However, for the on-orbit collection case, the velocities employed for the ejection are much larger than the ones induced by natural electrostatic forces. In the case of on-ground collection the ejection takes place in the dark side of the asteroid close to sunrise, where none or little electrostatic charge would build up, so it is also unlikely to be a useful mechanism if no artificial charging and electric fields are generated. Nonetheless, natural electrostatic forces need to be taken into account in this case, as the ejection velocities are low enough that very small particles may remain levitating after ejection, and collection takes place on the Sun-lit side where electric charges are building up. In addition, due to the different magnetic and electrical properties of the regolith components (ferromagnetic for the iron-nickel particles, paramagnetic for pyroxenes, non-magnetic for plagioclases) the final portion of their trajectories near the collection point will be influenced differently by the building up electrical charge on the surface, depending on the accumulated charge for each particle and its size.

Of greater concern are the cohesive/adhesive inter-particle forces in contact, which in theory significantly contribute to the internal strength of rubble piles [42]. Due to cleaner surfaces in vacuum and thus better contact between particles, cohesion can be up to 5 orders of magnitude larger than SRP at 1 AU for mm sized particles [6]. Two or more regolith grains in contact after being lifted from the surface are likely to remain together. Depending on the varying packing efficiency or porosity of bundles of grains held together by cohesion [43], the effective area-to-mass ratio may vary significantly, and with it the SRP force.

The main limitation on the model and analysis presented with regards to cohesion is therefore the assumption that each dust grain is treated as an independent particle, neglecting the formation of aggregates or clusters of particles due to cohesion or adhesion. Particle aggregates will have varying area-to-mass ratio depending on the shape and structure of the aggregate, and the orientation with respect to the incoming radiation. For example, for the simplest case of two particles (see Figure 15), the effective area-to-mass ratio, and corresponding lightness number $\beta^*$, can vary between 0.5 and 1 times the original $\beta$ of a single particle due to the orientation of the particles. In a more complex three-dimensional aggregate of a larger number of particles, the effective lightness number decreases with the equivalent radius of the cluster $r_{eq}$, and increases with porosity $\Phi$ (for higher porosities the mass per volume decreases). The effect of shape and orientation, represented in the following equation by a factor $f$, can result in large variations that are difficult to predict.

$$\beta^* \sim f \frac{1}{1-\Phi} \frac{r}{r_{eq}} \beta \qquad (19)$$

Aggregates with very high porosities can be formed in low-gravity conditions [43], depending on the force ratio. The force for two particles of the same radius can be expressed as [6]:

$$F_c = \frac{A_H S^2}{48\Omega^2} \frac{r}{2} \qquad (20)$$

where $A_H$ is the Hamaker coefficient, $S$ the cleanliness ratio (which can be assumed close to 1 for vacuum conditions, 0.1 in 'dirtier' environments), and $\Omega$ is a constant equivalent to the diameter of a $O^{-2}$ ion (~1.32 $10^{-10}$ m).

The force ratio between cohesive and gravity forces, or Bond ratio, is of the order of $10^6$ for particles of 100 µm on a 100 m radius asteroid with a rotation period of 5 hours. This corresponds to porosities of the order of 0.9 according to Yang et al. [43], or if the packing efficiency definition of Valverde and Castellanos [44] is used, to values for this packing efficiency between 0.1 and 0.18. Due to these high porosities, the decrease in lightness number is not as dramatic as initially could be expected. Aggregates of equivalent radius up to 10 times a single grain can have area-to-mass ratios of the same order of independent particles. In general, the effective lightness number will be smaller than that of a single particle, and the graph in Figure 14b will be displaced towards the left. There will be further mixing due to the fact that aggregates may consist of particles of different materials. The fact that the attraction force between a small particle and a larger grain is always greater than the forces between equal sized small particles [42] could also inhibit the application to size separation.

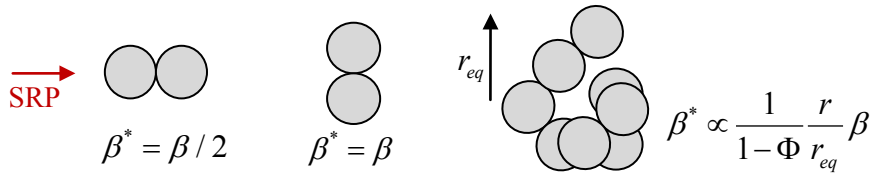

Figure 15. Effective lightness number for different configurations of aggregates.

In addition, complex models of particle cloud propagation may be required to take into account shading effects and collisions between ejected particles. Collisions may have the beneficial effect of evening out velocities at ejection. They may also break up clusters held together by cohesion if a joint receives an impact with sufficient kinetic energy, or, on the contrary, contribute to the build-up or restructuring of these aggregates, as shown in [45] and [46]. The complexity of these interactions is beyond the scope of this paper. Moreover, shading between particles may significantly reduce the SRP induced separation. Both effects are a function of the particle density in the cloud, which depends on the design of the ejection mechanism, and its mass flow rate.

Finally, bouncing and migration of particles after re-impact will affect the separation. If such a concentration method is implemented, special emphasis needs to be put on the design of the collectors to avoid undesired post-re-impact effects.

In general, cohesion represents the biggest drawback for the method proposed for finely ground particles, though it is unlikely that it will completely negate the segregation or concentration effect of solar radiation pressure as function of size or density. The strength of cohesive/adhesive forces is indeed an inherent problem for all separation techniques in vacuum and low-gravity environment, and other separation techniques that can be extrapolated from their terrestrial equivalent to the asteroid environment would suffer the same impediments.

**CONCLUSIONS**

The engineering of the solar radiation pressure (SRP) perturbation is a promising method for separation of dust grains around small bodies. Simplified models to describe the behavior of particles ejected at low velocities from an asteroid surface have been described and applied to a spherical rotating asteroid. The planar Hamiltonian approach and the phase space graphs introduced have proven to be useful tools to study and understand the behavior of dust, and they allow predicting the conditions to perform multiple revolutions or to re-impact as a function of the ejection site and the size of the particle.

A novel passive SRP separation method has been proposed, and possible variations of the collection strategy were discussed, both on-orbit and on-ground. The analysis suggests that the winnowing-like method with collection on ground can be an effective mechanism for material processing, while the on-orbit collection presents greater challenges. The method has the benefits of simplicity and a potential for high throughput, with possible applications for an initial pre-concentration of regolith sizes or materials prior to more complex processing methods. However, the efficacy of such method would greatly depend in the properties of the material, the conditions in each particular asteroid, and the effect of inter-particle forces that have not been taken into consideration in this paper. Also, surface mobility on asteroids in microgravity conditions represents a technological challenge for any collection strategy. Further analysis and demonstration of the concept through simulations with more complex models and/or microgravity and vacuum chamber laboratory tests would be of great value to assess the viability of SRP induced material separation.

## ACKNOWLEDGEMENTS

The work reported in this paper was supported by European Research Council grant 227571 (VISIONSPACE). The authors wish also to acknowledge the anonymous reviewers for their detailed and helpful comments and suggestions to improve the manuscript.